\documentclass{elsart}
\usepackage{graphics}
\usepackage{epsfig}
\usepackage{color}

\def\deg{{\rm deg}}

\def\qed{\hfill  \framebox(5,5){}}

\expandafter\chardef\csname pre amssym.def at\endcsname=\the\catcode`\@ \catcode`\@=11

\def\undefine#1{\let#1\undefined}
\def\newsymbol#1#2#3#4#5{\let\next@\relax
 \ifnum#2=\@ne\let\next@\msafam@\else
 \ifnum#2=\tw@\let\next@\msbfam@\fi\fi
 \mathchardef#1="#3\next@#4#5}
\def\mathhexbox@#1#2#3{\relax
 \ifmmode\mathpalette{}{\m@th\mathchar"#1#2#3}%
 \else\leavevmode\hbox{$\m@th\mathchar"#1#2#3$}\fi}
\def\hexnumber@#1{\ifcase#1 0\or 1\or 2\or 3\or 4\or 5\or 6\or 7\or 8\or
 9\or A\or B\or C\or D\or E\or F\fi}

\font\tenmsa=msam10 \font\sevenmsa=msam7 \font\fivemsa=msam5
\newfam\msafam
\textfont\msafam=\tenmsa \scriptfont\msafam=\sevenmsa \scriptscriptfont\msafam=\fivemsa \edef\msafam@{\hexnumber@\msafam}
\mathchardef\dabar@"0\msafam@39
\def\dashrightarrow{\mathrel{\dabar@\dabar@\mathchar"0\msafam@4B}}
\def\dashleftarrow{\mathrel{\mathchar"0\msafam@4C\dabar@\dabar@}}

        \font\tenmsb=msbm10
\font\sevenmsb=msbm7 \font\fivemsb=msbm5
\newfam\msbfam
\textfont\msbfam=\tenmsb \scriptfont\msbfam=\sevenmsb \scriptscriptfont\msbfam=\fivemsb \edef\msbfam@{\hexnumber@\msbfam}
\def\Bbb#1{\fam\msbfam\relax#1}

\newtheorem{theorem}{{\bf Theorem}}
\newtheorem{remark}{{\bf Remark}}
\newtheorem{corollary}[theorem]{{\bf Corollary}}
\newtheorem{proposition}[theorem]{{\bf Proposition}}
\newtheorem{lemma}[theorem]{{\bf Lemma}}

\begin{document}

\begin{frontmatter}



\title{Efficient Detection of Symmetries of Polynomially Parametrized Curves}


\author[a]{Juan Gerardo Alc\'azar\thanksref{proy}},
\ead{juange.alcazar@uah.es}

\address[a]{Departamento de Matem\'aticas, Universidad de Alcal\'a,
E-28871 Madrid, Spain}

\thanks[proy]{Supported by the Spanish `` Ministerio de
Ciencia e Innovacion" under the Project MTM2011-25816-C02-01.
Member of the Research Group {\sc asynacs} (Ref. {\sc ccee2011/r34}) }


\begin{abstract}
We present efficient algorithms for detecting central and mirror symmetry for the case of algebraic curves defined by means of polynomial parametrizations. The algorithms are based on the existence of a linear relationship between two proper polynomial parametrizations of the curve, which leads to a triangular polynomial system (with complex unknowns) that can be solved in a very fast way; in particular, curves parametrized by polynomials of serious degrees can be analyzed in a few seconds. In our analysis we provide a good number of theoretical results on symmetries of polynomial curves, algorithms for detecting rotation and mirror symmetry, and closed formulae to determine the symmetry center and the symmetry axis, when they exist. A complexity analysis of the algorithms is also given.
\end{abstract}
\end{frontmatter}

\section{Introduction}\label{section-introduction}

This paper deals with the problem of detecting the symmetries, if any, of a curve defined by means of a polynomial parametrization, i.e. a pair $(x(t),y(t))$ where both $x(t)$ and $y(t)$ are polynomials. The question of finding the symmetries of an algebraic curve has been previously investigated mainly because of its applications in pose estimation and patter recognition. Essentially, the problem is the following: a situation that is often studied in pattern recognition is how to choose, from a database of algebraic curves representing different objects, the one that best fits a given object, also represented by an algebraic equation; this question is addressed, among many others, in \cite{Huang} (where the database simulates different aircraft prototypes), \cite{Lei}, \cite{Sener}, \cite{TC00} (where the database corresponds to silhouettes of sea animals), \cite{Tasdizen} or \cite{Taubin1}. In order to do this one has to detect if the curve to be recognized is in fact the result of applying an affine transformation (typically translations, rotations, etc.) to some curve in the database; in turn, this implies to set the curve in a ``canonical position" that makes recognition possible. The problem of detecting symmetries comes along then as a means for setting properly the curve. In particular, this question has been addressed in \cite{Huang}, using splines, in \cite{Boutin}, \cite{Calabi}, \cite{Weiss} by means of differential invariants, in \cite{LR08}, \cite{LRTh}, \cite{TC00} using a complex representation of the implicit equation of the curve, or in \cite{Huang}, \cite{Suk1}, \cite{Suk2}, \cite{Taubin2} using moments.

It is also worth observing that in most of these papers it is assumed that the curve is given by means of its implicit equation. Exceptions to this are \cite{Huang}, where the curve is assumed to be a spline (i.e. a union of pieces of polynomial parametrizations), and \cite{Boutin}, \cite{Calabi} and other several papers on differential invariants, where the input is considered to be a parametrization without any restriction on its functional form. Furthermore, in almost all the papers on the question, the algorithm that is provided to detect symmetries is basically numerical, and therefore the output is approximate, and not exact. This is not a problem when the form to be recognized is, up to a certain extent, fuzzy, or when there are missing data, common situations in pattern recognition. But it is undesirable if the input is exact. Up to our knowledge, the exceptions to this are the papers on differential invariants (\cite{Boutin}, \cite{Calabi}, etc.) and \cite{LR08}, \cite{LRTh}. However, when applied to produce an exact output, the former can
only deal with parametrizations of low degree (in fact, the ultimate idea in the papers on differential invariants is to adapt the underlying theory to a numerical framework). As for the latter, i.e. \cite{LR08}, \cite{LRTh}, they provide deterministic algorithms for implicit algebraic curves $f(x,y)=0$, and yield exact answers in an efficient and elegant manner. However, if the curve is defined in parametric form, it is usual to try to keep in the same representation. In fact, even if we can convert from parametric to implicit form, this process is more complicated as the degree increases, and may produce non-sparse implicit equations with huge coefficients that are not convenient at all.

In this paper, we address the problem from a different perspective. On one hand, we focus on curves defined by means of polynomial parametrizations. Even if this seems too restrictive, polynomial parametrizations are widely used, for instance, in Solid Modeling and Computer Aided Geometric Design, and stay at the core of the notion of spline curve, commonly used in many applications. Furthermore, the algorithms that we provide can be directly conducted from the parametrization, and therefore do not require to compute the implicit equation of the curve (costly or even impossible for high degrees). Finally, we assume that our input is exact (the coefficients of the parametrization are required to be real numbers, not necessarily rational) and we provide also exact, i.e. deterministic, algorithms for checking whether the curve is symmetric or not, and for determining the elements of the symmetry, in the affirmative case.

The algorithms that we provide are very efficient and are able of analyzing curves with serious degrees in just a few seconds. The idea behind these algorithms comes from Real Algebraic Geometry, and exploits an algebraic property linking two ``good" (in a sense that is introduced in Section \ref{central}) parametrizations of a curve. In fact, in our case this relationship turns out to be linear. By applying this property, and writing the parametrization in complex form (which is really useful here) we are led to polynomial systems (with complex unknowns) but which are triangular. Hence, they can solved in a fast and efficient way. In fact, our results lead to closed formulae for the symmetry center and the symmetry axis, in the cases when they exist.

So, in the sequel we will analyze rotation and central symmetry (i.e. symmetry with respect to a point) in Section \ref{central}, and mirror symmetry (i.e. symmetry with respect to a line) in Section \ref{sec-mirror}. In Section \ref{float} we report on examples and timings showing the efficiency of our algorithms, and provide an analysis of the complexity. In Section \ref{conclusion} we provide some ideas on future work. Finally, in Appendix I (Section \ref{Appendix}) we list the examples mentioned in the paper, and provide information on them.

\section{Rotation and Central Symmetry} \label{central}

Along the paper we let ${\mathcal C}$ be an algebraic curve admitting a real polynomial parametrization $\varphi(t)=(x(t),y(t))$, i.e. $x(t),y(t)\in {\Bbb R}[t]$ (we will summarize this by saying that ${\mathcal C}$ is a {\it real polynomial curve}). Furthermore, we will also require $\varphi(t)$ to be {\it proper}, i.e. injective for almost all values of $t$. For instance, $(t,t^2)$ is a proper parametrization of a parabola, while $(t^2,t^4)$ is not (notice that whenever $t$ moves over the complex numbers, the latter parametrizes the parabola as well, and almost all the points of the parabola are generated by two different values of $t$). The properness of a parametrization is considered to be a good property, since non-proper parametrizations provide a ``redundant" representation of the curve, that might cause problems when plotting the curve, intersecting it with other curves, or just determining notable points of it.
From Theorem 6.11 in \cite{SWPD}, it is guaranteed that a polynomial curve can always be properly and polynomially reparametrized; more than that, we can ensure that in fact one can always find a proper polynomial parametrization over the reals. Indeed, if $\varphi(t)=(x(t),y(t))$ is polynomial but it is not proper, from the algorithm in page 193 of \cite{SWPD} we have that it can be properly reparametrized without extending the ground field; furthermore, if this reparametrization is not polynomial, then the algorithm in page 199 of \cite{SWPD} leads to a polynomial reparametrization, again without extending the ground field. So, in the rest of the paper, without loss of generality we will assume that $\varphi(t)$ is proper.

\subsection{Rotation Symmetry and Prohibitions.} \label{subsec-prb}

We say that ${\mathcal C}$ has {\sf rotation symmetry} iff there exists a point $P_0\in {\Bbb R}^2$ and an angle $\phi\in (0,2\pi)$ such that the rotation of center $P_0$ and angle $\phi$ leaves ${\mathcal C}$ invariant. In the special case when $\phi=\pi$, we say that ${\mathcal C}$ has {\sf central symmetry}, i.e. that ${\mathcal C}$ is symmetric with respect to $P_0$ (the center of symmetry). The following theorem proves that in our case, the only form of rotation symmetry that a polynomial curve can exhibit is central symmetry. Here we use the notion of {\it infinite branch} of a parametric curve: by this notion, we mean intuitively a part of ${\mathcal C}$ that goes to infinity; for example, the parabola $(t,t^2)$ has, according to this terminology, two infinite branches, one for $t=\infty$ and another one for $t=-\infty$. In fact, every polynomial curve has only two infinite branches; however, rational curves may have more than two (because there may be $t$-values where the denominator of either $x(t)$ or $y(t)$ vanishes).

\begin{theorem} \label{only-central}
Real polynomial curves cannot have any other form of rotation symmetry other than central symmetry.
\end{theorem}

{\bf Proof.} if ${\mathcal C}$ presents rotation symmetry other than central symmetry then it consists of the union of $m$ different copies, with $m\geq 3$. Since each copy provides at least one infinite branch, we get at least $m$ infinite branches. But $m\geq 3$, when the
number of infinite branches of ${\mathcal C}$ is 2. \qed

Furthermore, the following result provides a necessary condition for ${\mathcal C}$ to have central symmetry.

\begin{theorem} \label{sufficient}
Let ${\mathcal C}$ be a real polynomial curve, and let $\varphi(t)=(x(t),y(t))$ be a proper, real, polynomial parametrization of ${\mathcal C}$. If ${\mathcal C}$ has central symmetry, then $\mbox{lim}_{t\to \infty}x(t)$ and $\mbox{lim}_{t\to -\infty}x(t)$, $\mbox{lim}_{t\to \infty}y(t)$ and $\mbox{lim}_{t\to -\infty}y(t)$, exhibit different signs. In particular, $\deg_t(x(t))$ and $\deg_t(y(t))$ must be both odd.
\end{theorem}

{\bf Proof.} If ${\mathcal C}$ has central symmetry then it must exhibit two different infinite branches, that go to infinity in opposite quadrants. \qed

\subsection{Detecting Central Symmetry} \label{subsec-detect}

According to the preceding subsection, central symmetry is the only form of rotation symmetry that ${\mathcal C}$ can exhibit. So, let us address here the problem of detecting this kind of symmetry in an efficient way. For this purpose, previously we need the following result, that is crucial for us.

\begin{lemma} \label{z-prop}
Let $\varphi_1(t)$ and $\varphi_2(t)$ be two polynomial and proper parametrizations of the same curve. Then there exists a linear function $L(t)=\alpha t +\beta$, with $\alpha,\beta\in {\Bbb R}$, such that $\varphi_1(t)=\varphi_2(\alpha t + \beta)$ and $\alpha\neq 0$.
\end{lemma}

{\bf Proof.} By Lemma 4.17 in \cite{SWPD}, $\varphi_1(t)$ and $\varphi_2(t)$ can be obtained from the other by means of a rational reparametrization $\xi(t)=\displaystyle{\frac{\alpha t + \beta}{\gamma t + \delta}}$, where $\alpha,\beta,\gamma,\delta\in {\Bbb C}$ and $\alpha \delta-\beta  \gamma\neq 0$. Since both $\varphi_1(t)$ and $\varphi_2(t)$ are polynomial, from Lemma 6.8 in \cite{SWPD} we deduce that in fact $\xi(t)$ must be {\it linear}. Furthermore, let us see that $\alpha,\beta\in {\Bbb R}$. Indeed, if $\alpha,\beta\in {\Bbb C}$ then
  since $\varphi_1(t)=\varphi_2(\alpha t + \beta)$ and $\varphi_1(t)$ is a real polynomial parametrization, we have that $\varphi_2(\overline{\alpha}+t\overline{\beta})=\varphi_2(\alpha t + \beta)$ for all $t\in {\Bbb R}$. But then $\varphi_2(t)$ cannot be proper, because for every (real or complex) value of $t$, it holds that $\alpha t+\beta$ and $\overline{\alpha}t+\overline{\beta}$ generate the same point. Since $\varphi_2(t)$ is proper by hypothesis, then $\alpha,\beta\in {\Bbb R}$. Finally, since in this case $\delta=1$ and $\gamma=0$, then $\alpha \delta-\beta \gamma=\alpha \neq 0$. \qed

Along the paper, it will be more convenient to
write $\varphi(t)$ in complex form as $z(t)=x(t)+{\bf i}\cdot y(t)$. Now ${\mathcal C}$ has central symmetry, with center of symmetry $z_0$ ($z_0$ corresponds to a point $P_0=(x_0,y_0)$, written in complex form), if and only if for every $t$-value, $z^{\star}(t)=-z(t)+2z_0$ is also a point of ${\mathcal C}$ (see Figure 1), i.e. if and only if $z^{\star}(t)$ is another  parametrization of ${\mathcal C}$ (in complex form). Now one may check that $z^{\star}(t)$ also corresponds to a real polynomial parametrization, different from $z(t)$; furthermore, since $z^{\star}(t)$ is the result of composing $\varphi(t)$ with a bijective planar transformation (the symmetry with respect to $z_0$), then it is also a proper parametrization. So, ${\mathcal C}$ has central symmetry if and only if $z(t)$ and $z^{\star}(t)$ are two {\it proper} parametrizations of the same curve. These observations, together with Lemma \ref{z-prop}, give rise to the following theorem.

\begin{figure}[ht]
\begin{center}
\centerline{$\begin{array}{c}
\includegraphics[width=7cm,height=6cm]{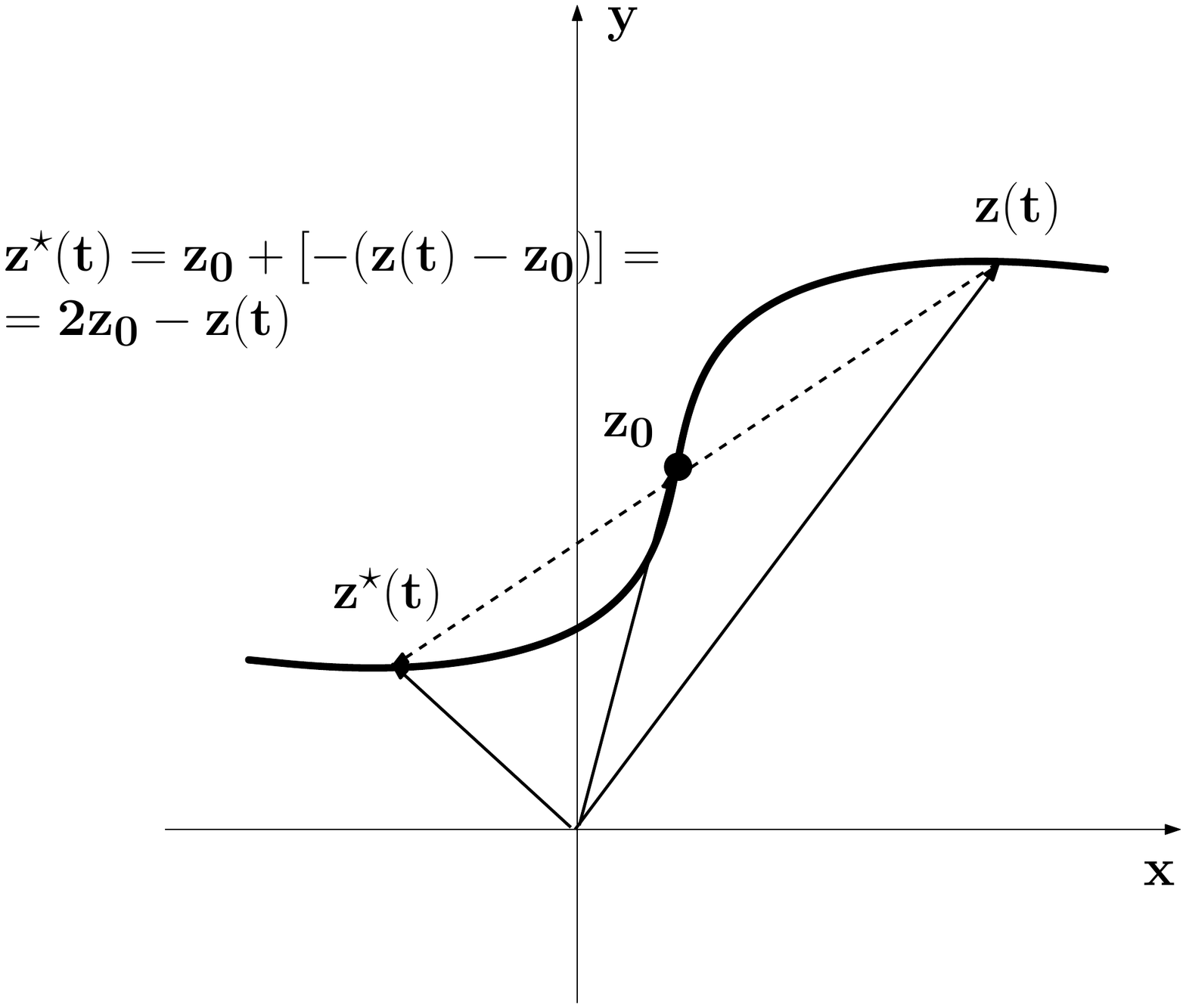}
\end{array}$}
\end{center}
\caption{Central Symmetry}
\end{figure}

\begin{theorem} \label{th-main}
The curve ${\mathcal C}$ has central symmetry if and only if there exist $z_0\in {\Bbb C}$ and a linear transformation $L(t)=\alpha t + \beta$, with $\alpha,\beta\in {\Bbb R}$, $\alpha\neq 0$, such that $z^{\star}(t)=z(L(t))$.
\end{theorem}



Now writing $z(t)=c_nt^n+c_{n-1}t^{n-1}+\cdots + c_0$, where $c_i\in {\Bbb C}$ for $i=0,1,\ldots,n$, and $L(t)=\alpha t +\beta$, with $\alpha,\beta\in {\Bbb R}$, the statement in Theorem \ref{th-main} leads to the following $(n+1) \times (n+1)$ polynomial system, that we denote as ${\mathcal S}$:
\begin{equation} \label{eq-central}
\hspace*{-1in}\begin{array}{rcccccccccc}
-(c_0-z_0) + z_0 & = & c_0 & +& c_1\beta & + & c_2\beta^2 & + \cdots + & c_{n-1}\beta^{n-1} & + & c_n\beta^n \\
-c_1 & = &  & & c_1\alpha & + & c_2\cdot 2\alpha\beta & + \cdots +  & c_{n-1}(n-1)\alpha \beta^{n-2} & + & c_n n \alpha \beta^{n-1}\\
-c_2 & = &  & &  & & c_2\alpha^2 & + \cdots +  & c_{n-1}\frac{(n-1)(n-2)}{2}\alpha^2 \beta^{n-3} & + & c_n \frac{n(n-1)}{2} \alpha^2 \beta^{n-2}\\
  & \vdots &  & &   &   &   &   & \vdots  &   & \vdots \\
-c_{n-1} & = & &  &  & &  &   &  c_{n-1}\alpha^{n-1} & + & c_n n\alpha^{n-1}\beta \\
-c_{n} & = & &  &  & &  &   &    &   & c_n \alpha^{n}
\end{array}
\end{equation}
Since by Theorem \ref{sufficient} the value of $n$ must be odd, and the case $n=1$ is trivial (${\mathcal C}$ is a line), in the rest of the section we will assume that $n\geq 3$. Now we will refer to the above equations as $(0),(1),\ldots,(n)$, respectively (i.e. the equation $(k)$ of ${\mathcal S}$ is the one containing $\alpha^k$). Observe that Theorem \ref{th-main} is equivalent to the existence of a solution $(z_0,\alpha,\beta)$, with $\alpha,\beta\in {\Bbb R}$, of the system ${\mathcal S}$. This system has a triangular structure that can be exploited for solving the problem in an efficient way: first of all, from the last equation we have $\alpha^n=-1$, with $n$ odd. So, $\alpha=-1$. Additionally, from the equation $(n-1)$ we deduce that $\beta=\frac{-2c_{n-1}}{n\cdot c_n}$, which must be real in the case when
${\mathcal C}$ is symmetric. So, the only possible value for $(\alpha,\beta)$ is $(-1,\frac{-2c_{n-1}}{n\cdot c_n})$. Finally, since the coefficient of $z_0$ in the equation (0) is $-2$, whenever $(-1,\frac{-2c_{n-1}}{n\cdot c_n})$ fulfill the equations $(1),\ldots,(n)$ we can solve for $z_0$ in the equation (0), and obtain the symmetry center. We summarize these ideas in the following theorem, where a closed expression for the symmetry center (in the case when it exists) is provided.




\begin{theorem} \label{main-rot}
The curve ${\mathcal C}$ has central symmetry if and only if $\beta=\frac{-2c_{n-1}}{n\cdot c_n}\in {\Bbb R}$, and $\alpha=-1,\beta=\frac{-2c_{n-1}}{n\cdot c_n}$ satisfy the equations $(1),\ldots,(n)$ of the system ${\mathcal S}$. Moreover, in that case the center of symmetry is the point $z_0$ (in complex form) given by
\begin{equation} \label{center}
z_0=c_0+\displaystyle{\frac{1}{2}(c_1\beta+c_2\beta^2+\cdots +c_n\beta^n)}
\end{equation}
\end{theorem}


\begin{remark} \label{arecomplex}
From Theorem 5.3 in \cite{LRTh}, it follows that ${\mathcal C}$ can have at most just one center of symmetry. In our case, the existence of at most one symmetry center follows also from the above theorem.
\end{remark}

This theorem provides in a natural way an algorithm for checking central symmetry, and for determining the center of symmetry in the affirmative case. This is illustrated by the following example.

{\bf Example 1.} {\it Consider the curve ${\mathcal C}$ parametrized by $(x(t),y(t))$, where
\[
\left\{\begin{array}{rcl}
x(t) & = & 2+2(2t+1)^{23}-(2t+1)^{13}+2(2t+1)^{11}+2(2t+1)^5-(2t+1)^3+2t\\
y(t) & = & -2(2t+1)^{23}+(2t+1)^{13}-2(2t+1)^{11}+2(2t+1)^5-(2t+1)^3+2t
\end{array}\right.
\]
One can check that the curve is proper by applying Theorem 4.30 in \cite{SWPD}. It is perhaps instructive to mention that we obtained this curve by considering first $\psi(t)=(2t^{23}-t^{13}+2t^{11},2t^{5}-t^3+t)$, which is obviously symmetric with respect to the origin, then substituting $t:=2t+1$, and finally applying the change of coordinates $\{X=1+x+y,Y=-1-x+y\}$, which is the composition of a rotation and a translation of vector $(-1,1)$; as a consequence, we reach a curve which is symmetric with respect to the point $(-1,1)$. Now one can check that the implicit equation $f(x,y)$ of ${\mathcal C}$ has degree 23 and that the infinity norm of this implicit equation is close to $2^{600}$. The size of the implicit equation suggests that it would be difficult to work with it in implicit form. Using our method and Maple 14, running in a computer with 8 Gb of RAM and a CPU revving up to 2 GHz., it takes 0.374 seconds to build the system ${\mathcal S}$ and check that $\alpha=-1$, $\beta=-1$, $z_0=1-i$ is a solution. So, we recover $(1,-1)$ as the symmetry center of the curve. }

In the case when $c_{n-1}=0$ the analysis of the system ${\mathcal S}$ is easier. Indeed, in this case from the equation $(n-1)$ we get $\beta=0$. Hence, we have the following result.

\begin{theorem}\label{spec-case}
Let ${\mathcal C}$ be a real polynomial curve and let $z(t)=c_nt^n+c_{n-1}t^{n-1}+\cdots +c_1t+c_0$ be a polynomial, proper, real parametrization of ${\mathcal C}$ in complex form. If $c_{n-1}=0$, then the following statements hold:
\begin{itemize}
\item [(i)] If $n=3$ (i.e. ${\mathcal C}$ is a cubic) then ${\mathcal C}$ is symmetric with respect to the point $c_0$.
\item [(ii)] If $n\geq 4$ then ${\mathcal C}$ has central symmetry if and only if $c_k=0$ when $k\in \{1,\ldots,n-2\}$ is even. Moreover, in that case the symmetry center is $c_0$.
\end{itemize}
\end{theorem}

{\bf Proof.} In the case of (i), one can check that $\alpha=-1,\beta=0$ fulfill the equations (1), (2), (3) of the system ${\mathcal S}$. Then plugging $\alpha=-1,\beta=0$ into (0) we get $z_0=c_0$. In the case of (ii), when we substitute $\beta=0$ in the equations $(1),\ldots,(n-2)$ we get $c_k\cdot(-1)^k+c_k=0$ for $k=1,\ldots,n-2$. When $k$ is odd this equality is clearly fulfilled, and when $k$ is even the equality holds iff $c_k=0$. \qed

\section{Mirror Symmetry} \label{sec-mirror}

 \subsection{Detecting Mirror Symmetry} \label{subsec-detetc}

 Along this section, we write the {\it conjugate} of $z(t)$ as $\overline{z(t)}=x(t)-{\bf i}\cdot y(t)$. Now we say that ${\mathcal C}$ has {\sf mirror symmetry} if there exists an axis ${\mathcal L}$ (called the {\it symmetry axis}) such that ${\mathcal C}$ is symmetric with respect to ${\mathcal L}$. The curve ${\mathcal C}$ has this type of symmetry if and only if there exists a movement $M$, composition of a translation and a rotation, such that $M({\mathcal C})$ is symmetric with respect to the $x$-axis; that is to say, if and only there exist $z_0\in {\Bbb C}$, and $\phi\in [0,2\pi)$, such that $\tilde{z}(t)=(z(t)-z_0)\cdot e^{{\bf i}\cdot \phi}$ is, in complex form, the parametrization of a curve $\tilde{\mathcal C}$, symmetric with respect to the $x$-axis. So, ${\mathcal C}$ exhibits mirror symmetry if and only if $w(t)=(z(t)-z_0)\cdot e^{{\bf i}\cdot \phi}$ and $\overline{w}(t)=\overline{(z(t)-z_0)\cdot e^{{\bf i}\cdot \phi}}$ both parametrize the same curve $\tilde{\mathcal C}$. If $z(t)$ corresponds to a proper parametrization then $w(t)$ and $\overline{w}(t)$ define also proper parametrizations of the corresponding curves. Hence, taking also into account Lemma \ref{z-prop}, the following theorem, analogous to Theorem \ref{th-main} in Section \ref{central}, holds.

\begin{figure}[ht]
\begin{center}
\centerline{$\begin{array}{c}
\includegraphics[width=7cm,height=6cm]{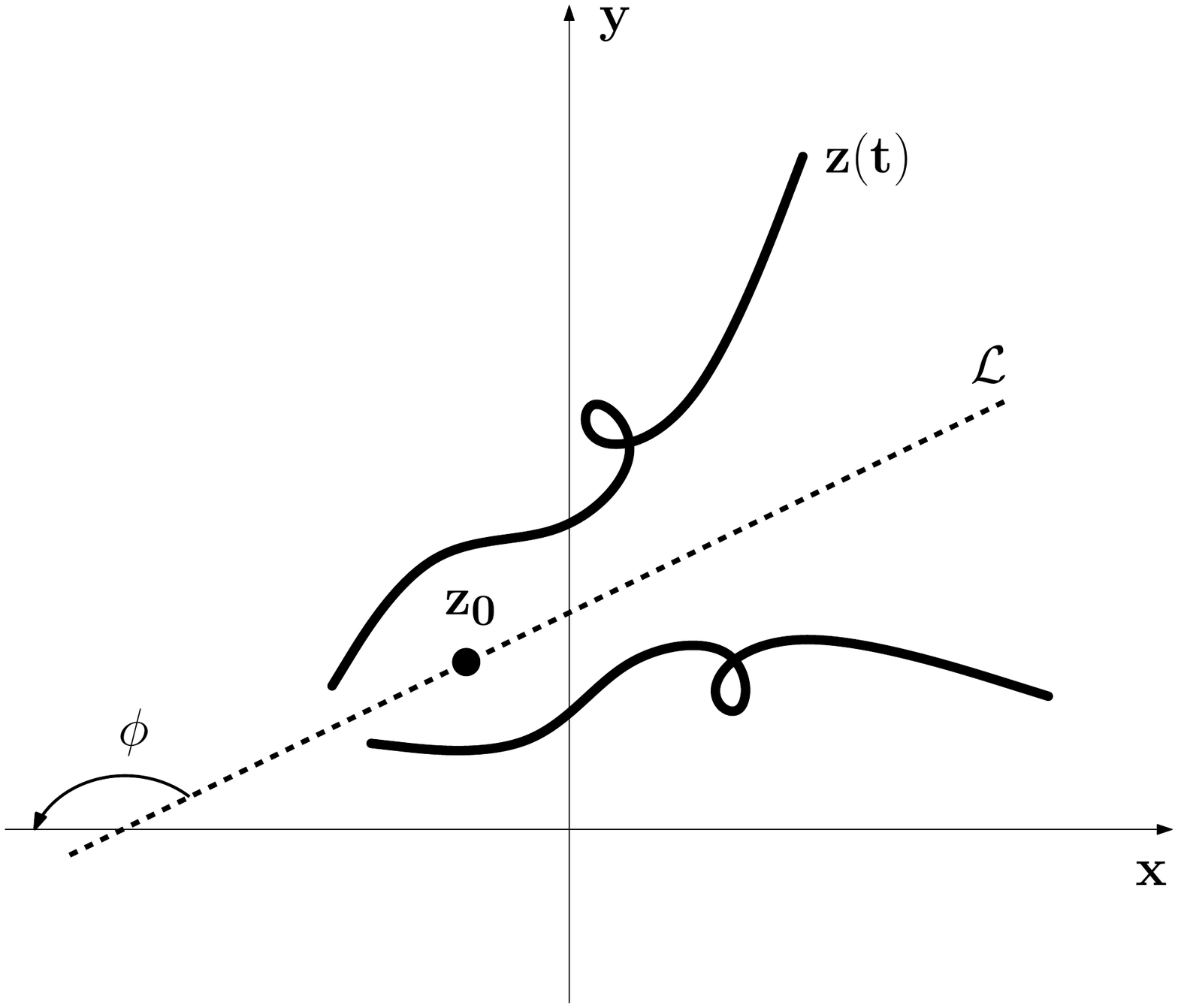}
\end{array}$}
\end{center}
\caption{Mirror Symmetry}
\end{figure}

\begin{theorem} \label{th-main-2}
The curve ${\mathcal C}$ exhibits mirror symmetry if and only if there exist a line ${\mathcal L}$ and
$\alpha,\beta\in {\Bbb R}$, $\alpha\neq 0$, such that for any $z_0\in {\mathcal L}$, it holds that $\overline{w}(t)=w(\alpha t +\beta)$.
\end{theorem}

Let us see how to take advantage of this theorem for efficiently checking, and computing in the affirmative case, the existence of mirror symmetry. We need first the following lemma.

\begin{lemma} \label{one-most}
A polynomial curve cannot have more than one symmetry axis.
\end{lemma}

{\bf Proof.} If ${\mathcal C}$ has $m$ symmetry axes then by symmetry each infinite branch gives us $m$ infinite branches more. However, since ${\mathcal C}$ is polynomial then it has only two infinite branches, and therefore $m=1$. \qed

We will also need the following result, where we use the notation and terminology of Theorem \ref{th-main-2}.

\begin{lemma} \label{valueof}
If $\overline{w}(t)=w(\alpha t +\beta)$ for $\alpha,\beta\in {\Bbb R}$, then $\alpha=-1$.
\end{lemma}

{\bf Proof.} If $\overline{w}(t)=w(\alpha t +\beta)$ for $\alpha,\beta\in {\Bbb R}$, then by Theorem \ref{th-main-2} the curve ${\mathcal C}$ is symmetric with respect to a certain axis, and $w(t)$ parametrizes a curve $\tilde{\mathcal C}$ which is symmetric with respect to the $x$-axis. Let us write $w(t)=(a_pt^p+a_{p-1}t^{p-1}+\cdots, b_qt^q+b_{q-1}t^{q-1}+\cdots)$. Because of the symmetry of $\tilde{\mathcal C}$ with respect to the $x$-axis, we have that $q$ must be odd (so that the sign of $y(t)$ when $t$ goes to $\infty$ or $-\infty$ is different). Now since $w(\alpha t+\beta)=\overline{w}(t)$, then \[b_q\alpha^q t^q+\cdots=-b_q t^q+\cdots\]Hence, $\alpha^q=-1$. Since $q$ is odd, $b_q\neq 0$ and $\alpha$ is real, then $\alpha=-1$. \qed

Now taking into account Theorem \ref{th-main-2} and the expression of $w(t)$, we get the following system ${\mathcal W}$:
\begin{equation}
\hspace*{-1in}\begin{array}{rcccccccccc}
(\overline{c}_0-\overline{z}_0)\cdot e^{-{\bf i}\cdot 2\phi} & = & (c_0-z_0) & +& c_1\beta & + & c_2\beta^2 & + \cdots + & c_{n-1}\beta^{n-1} & + & c_n\beta^n\\
\overline{c}_1\cdot e^{-{\bf i}\cdot 2\phi} & = &  & & c_1\alpha & + & c_2\cdot 2\alpha\beta & + \cdots +  & c_{n-1}(n-1)\alpha \beta^{n-2} & + & c_n n \alpha \beta^{n-1}\\
\overline{c}_2\cdot e^{-{\bf i}\cdot 2\phi} & = &  & &  & & c_2\alpha^2 & + \cdots +  & c_{n-1}\frac{(n-1)(n-2)}{2}\alpha^2 \beta^{n-3} & + & c_n \frac{n(n-1)}{2} \alpha^2 \beta^{n-2}\\
  & \vdots &  & &   &   &   &   & \vdots  &   & \vdots \\
\overline{c}_{n-1}\cdot e^{-{\bf i}\cdot 2\phi} & = & &  &  & &  &   &  c_{n-1}\alpha^{n-1} & + & c_n n\alpha^{n-1}\beta \\
\overline{c}_{n}\cdot e^{-{\bf i}\cdot 2\phi} & = & &  &  & &  &   &    &   & c_n \alpha^{n}
\end{array}
\end{equation}
We denote the above equations as $(0),(1), \ldots, (n)$, respectively (i.e. the equation $(k)$ of ${\mathcal W}$ is the one containing the power $\alpha^k$). Notice that the equation (0) is the only one containing $z_0$, and is linear in $z_0$. Furthermore, from Theorem \ref{th-main-2} one can see that the existence of mirror symmetry is equivalent to the existence of $\phi\in [0,2\pi)$, $\alpha,\beta\in {\Bbb R}$ such that: (i) these values fulfill the equations $(1),\ldots,(n)$; (ii) for these values, the equation (0) corresponds to a real line (which would be the symmetry axis ${\mathcal L}$). Now by Lemma \ref{valueof}, if ${\mathcal W}$ is consistent then $\alpha=-1$; also, by combining the last two equations we get that
\begin{equation} \label{b-1}
\beta=\xi(\alpha)=\displaystyle{\frac{c_n \bar{c}_{n-1}\alpha-\bar{c}_n c_{n-1}}{n|c_n|^2}}
\end{equation}
And since $\alpha=-1$, we get
\begin{equation} \label{b-2}
\beta=\xi(-1)=\displaystyle{-\frac{c_n \bar{c}_{n-1}+\bar{c}_n c_{n-1}}{n|c_n|^2}}=\displaystyle{\frac{-2\cdot \mbox{Re}(c_n\cdot \bar{c}_{n-1})}{n|c_n|^2}}
\end{equation}
Observe that $\beta\in {\Bbb R}$. Now in order to verify whether the equation (0) of ${\mathcal W}$ corresponds to a real line, the following results are useful. The first one can be proven in a straightforward way.

\begin{lemma} \label{lem-line}
Any real line $Ax+By+C=0$ is transformed by means of the complex change $x=\displaystyle{\frac{z+\bar{z}}{2}}$, $y=\displaystyle{\frac{z-\bar{z}}{2i}}$ into $\overline{\gamma}z+\gamma \bar{z}+C=0$, where $\gamma=\frac{A}{2}+i\frac{B}{2}$. Conversely, any equality $\overline{\gamma}z+\gamma \bar{z}+C=0$ with $C$ real corresponds to a real line.
\end{lemma}

The above lemma is used for proving the following result. Here, we denote $Q(\beta)=c_1\beta+c_2\beta^2+\cdots+c_n\beta^n$, and
\begin{equation} \label{b-3}
Q^{\star}(\beta)=Q(\beta)\cdot \displaystyle{\frac{\bar{c}_n}{|c_n|}\cdot \sqrt{(-1)^{n+1}}}
\end{equation}

\begin{proposition} \label{real-line}
The equation $(0)$ of the system ${\mathcal W}$ corresponds to a real line if and only if $Q^{\star}(\beta)\in {\Bbb R}$.
\end{proposition}

{\bf Proof.} After substituting $-e^{-{\bf i}2\phi}=e^{{\bf i}(\pi-2\phi)}=e^{{\bf i}2(\frac{\pi}{2}-\phi)}$ in the equation (0) of ${\mathcal W}$, and dividing the whole equation by $e^{{\bf i}(\frac{\pi}{2}-\phi)}$, we get \[e^{-{\bf i}(\pi/2-\phi)}z_0+e^{{\bf i}(\frac{\pi}{2}-\phi)}\bar{z}_0=\bar{c}_0\cdot e^{{\bf i}(\frac{\pi}{2}-\phi)}+c_0\cdot e^{-{\bf i}(\frac{\pi}{2}-\phi)}+Q(\beta)\cdot e^{-{\bf i}(\frac{\pi}{2}-\phi)}\]Now from
Lemma \ref{lem-line}, and taking into account that the conjugate of $e^{-{\bf i}(\pi/2-\phi)}$ is $e^{{\bf i}(\pi/2-\phi)}$, it holds that the above expression corresponds to a real line iff the right hand-side of the equation is a real number. However, one may see that
\[\bar{c}_0\cdot e^{{\bf i}(\frac{\pi}{2}-\phi)}+c_0\cdot e^{-{\bf i}(\frac{\pi}{2}-\phi)}\]is already real, because it is the sum of a complex number, and its conjugate. So, we get a real line iff $Q(\beta)\cdot e^{-{\bf i}(\frac{\pi}{2}-\phi)}$ is a real number. Finally, using that $-e^{-{\bf i}2\phi}=e^{{\bf i}2(\frac{\pi}{2}-\phi)}$ and the equation $(n)$, one has that
\[e^{-{\bf i}(\frac{\pi}{2}-\phi)}=\displaystyle{\frac{1}{e^{{\bf i}(\frac{\pi}{2}-\phi)}}=\frac{1}{\sqrt{-e^{-{{\bf i}2\phi}}}}=\frac{1}{\sqrt{\frac{-c_n (-1)^n}{\bar{c}_n}}}=\sqrt{\frac{-\bar{c}_n}{c_n(-1)^n}}}\]Multiplying and dividing by $\bar{c}_n$ into the square root, the result follows easily. \qed


So, we can put together all these observations to characterize the existence of mirror symmetry. For this purpose, it is preferable to consider a new system ${\mathcal W}'$, which is obtained from ${\mathcal W}$ by substituting
\begin{equation} \label{ei}
e^{-{\bf i}\cdot 2\phi}=\frac{c_n}{\overline{c}_n}\alpha^n
\end{equation}
(from the last equation of ${\mathcal W}$) into the $n$ first equations, and dividing out the equation $(k)$ by $\alpha^k$:
\begin{equation}
\hspace*{-1.1in}\begin{array}{rcccccccccc}
(\overline{c}_0-\overline{z}_0)c_n\cdot \alpha^n & = & \overline{c}_n\cdot [(c_0 -z_0)& +& c_1\beta & + & c_2\beta^2 & + \cdots + & c_{n-1}\beta^{n-1} & + & c_n\beta^n]\\
\overline{c}_1 c_n\cdot \alpha^{n-1} & = &  & & \overline{c}_n\cdot [ c_1 & + & 2c_2 \beta & + \cdots +  & c_{n-1}(n-1) \beta^{n-2} & + & c_n n  \beta^{n-1}]\\
\overline{c}_2 c_n\cdot\alpha^{n-2} & = &  & &  & & \overline{c}_n\cdot [ c_2  & + \cdots +  & c_{n-1}\frac{(n-1)(n-2)}{2}  \beta^{n-3} & + & c_n \frac{n(n-1)}{2}   \beta^{n-2}]\\
  & \vdots &  & &   &   &   &   & \vdots  &   & \vdots \\
\overline{c}_{n-1}c_n\cdot\alpha & = & &  &  & &  &   &  \overline{c}_n\cdot [ c_{n-1}   & + & c_n n \beta]
\end{array}
\end{equation}
We denote these equations as $[0], [1], \ldots , [n-1]$, respectively. Furthermore, we write the equation [0] as $r(z_0,\bar{z}_0,\beta)=0$. Then the following result, which can be deduced from the results in this section, holds.

\begin{theorem} \label{main-mirr}
The curve ${\mathcal C}$ has mirror symmetry if and only if the following two conditions hold: (1) $\alpha=-1$, $\beta=\xi(-1)$ fulfill the equations $[1],\ldots,[n-1]$; (2) $Q^{\star}(\beta)\in {\Bbb R}$. Moreover, in that case the symmetry axis of ${\mathcal C}$ (in complex form) is $r(z,\bar{z},\beta)$, i.e.
\begin{equation} \label{axis}
\overline{c}_nz-c_n(-1)^n\overline{z}+\bar{c}_0c_n(-1)^n-\bar{c}_nc_0-Q(\beta)=0
\end{equation}
\end{theorem}



Theorem \ref{main-mirr} provides also the following result on the existence of symmetries of ${\mathcal C}$.

\begin{proposition} \label{some-prohib}
The following statements are true:
\begin{itemize}
\item [(1)] If $(-1)^nc_n-\bar{c}_n=0$, and ${\mathcal C}$ has a symmetry axis, then it is parallel to the $x$-axis.
\item [(2)] If $(-1)^nc_n+\bar{c}_n=0$, and ${\mathcal C}$ has a symmetry axis, then it is parallel to the $y$-axis.
\end{itemize}
\end{proposition}

{\bf Proof.} Let us see (1). For this purpose, assume that ${\mathcal L}$ is a symmetry axis for ${\mathcal C}$. The intersection of ${\mathcal L}$ with the $x$-axis can be found by substituting $\bar{z}_0=z_0$ in the closed expression for $z_0$ provided in the statement of Theorem \ref{main-mirr}. This substitution yields $(\bar{c}_n-(-1)^nc_n)z_0+\cdots=0$ (where $\cdots$ comprises terms where $z_0$ is not present), and hence when $\bar{c}_n-(-1)^nc_n=0$ no value for $z_0$ can be found. So, (1) follows. The statement (2) is proven in a similar way but imposing $\bar{z}_0=-z_0$, instead of $\bar{z}_0=z_0$. \qed

The above proposition can be used for proving the following result on the existence of symmetries. Here we denote $\deg_t(x(t))=r$, $\deg_t(y(t))=s$.

\begin{theorem} \label{conclu-some-pro}
The following statements hold:
\begin{itemize}
\item [(i)] If $r>s$ and $r$ is even, and ${\mathcal C}$ has a symmetry axis, then it is parallel to the $x$-axis.
\item [(ii)] If $r>s$ and $r$ is odd, and ${\mathcal C}$ has a symmetry axis, then it is parallel to the $y$-axis.
\item [(iii)] If $r<s$ and $s$ is odd, and ${\mathcal C}$ has a symmetry axis, then it is parallel to the $x$-axis.
\item [(iv)] If $r<s$ and $s$ is even, and ${\mathcal C}$ has a symmetry axis, then it is parallel to the $y$-axis.
\end{itemize}
\end{theorem}

{\bf Proof.} If $r>s$ then $n=r$, and $c_n$ is real, i.e. $\bar{c}_n=c_n$. Now if $n$, i.e. $r$, is even we are in the hypotheses of statement (i), and the result follows from the statement (1) of Proposition \ref{some-prohib}; if $n$, i.e. $r$, is odd we we are in the hypotheses of statement (ii), which follows from the statement (2) of Proposition \ref{some-prohib}. Finally if $r<s$ then $n=s$, and $c_n$ is a purely imaginary number, i.e. $\bar{c}_n=-c_n$; then we argue as before. \qed

Theorem \ref{main-mirr} leads to an algorithm for computing the symmetry axis of ${\mathcal C}$ (if any). This algorithm is illustrated in the following example.

{\bf Example 2.} {\it Consider the curve ${\mathcal C}$ parametrized by $(x(t),y(t))$, where
\[
\left\{\begin{array}{rcl}
x(t) & = & (2t+1)^{20}+(2t+1)^{18}+(2t+1)^{10}+1+(2t+1)^{21}-3(2t+1)^5+(2t+1)^3\\
y(t) & = & -(2t+1)^{20}-(2t+1)^{18}-(2t+1)^{10}-1+(2t+1)^{21}-3(2t+1)^5+(2t+1)^3
\end{array}\right.
\]
This curve is proper, and has been constructed starting from a simpler curve, then applying a change of parameters and finally rotating the curve $\frac{\pi}{4}$ radians.
One can check that the implicit equation $f(x,y)$ of ${\mathcal C}$ has degree 21, that the infinity norm of this implicit equation is close to $2^{500}$. Using our method, it takes 1.388 seconds to
construct the system ${\mathcal W}'$ and check that:
\begin{itemize}
\item $\alpha=-1$, $\beta=\xi(-1)$, fulfill the equations $[1],\ldots,[n-1]$, and $Q^{\star}(\beta)=2\sqrt{2}\in {\Bbb R}$. So, from Theorem \ref{main-mirr}, ${\mathcal C}$ has one symmetry axis.
    \item The equation of the symmetry axis is obtained by substituting the above values for $\alpha,\beta$ in the equation [0], which yields
    $iz-\bar{z}=0$, i.e. $y=x$.
\end{itemize}}

\noindent
As in the preceding section, in the case $c_{n-1}=0$ we can provide a sharper result. In this case $\beta=0$ and therefore $Q^{\star}(\beta)=0$; so, the second condition in Theorem \ref{main-mirr} always holds. In addition to this, the equations of ${\mathcal W}'$ have the form $m_k(\alpha)=c_n\overline{c}_k\alpha^k-\overline{c}_n c_k=0$ for $k=1,\ldots,n-2$. So, when $\alpha=-1$ we have
\[m_k(\alpha)=\left\{\begin{array}{ccl}2{\bf i}\cdot \mbox{Im}(c_n\bar{c}_k) & \mbox{ if } & $k$ \mbox{ is even}\\
-2 \cdot \mbox{Re}(c_n\bar{c}_k) & \mbox{ if } & $k$ \mbox{ is odd}
\end{array}\right.\]

Hence, the next result, which essentially follows from the above expression, holds.

\begin{theorem} \label{spec-case-mirror}
Assume that $c_{n-1}=0$. Then ${\mathcal C}$ presents mirror symmetry if and only if one of for every $k=1,\ldots,n-2$, $c_n\cdot\bar{c}_k\in{\Bbb R}$ when $k$ is even, and $c_n\cdot\bar{c}_k$ is either a pure imaginary number (i.e. with null real part) or $0$ when $k$ is odd.
        \end{theorem}

\begin{remark} \label{indep} The above theorem shows that in fact the two conditions in Theorem \ref{main-mirr} are independent; in particular, in the above case it can happen that $Q^{\star}(\beta)\in {\Bbb R}$ but the curve does not have mirror symmetry.
\end{remark}

\subsection{Some More Prohibitions and Observations}\label{some-more}

If ${\mathcal C}$ exhibits mirror symmetry and ${\mathcal L}$ is a symmetry axis, then almost all elements in the family of lines $Ax+By+C=0$ normal to ${\mathcal L}$ must intersect ${\mathcal C}$ at an even number of points (the finitely many exceptions will be those lines that contain some point of ${\mathcal C}\cap {\mathcal L}$). We can exploit this for deducing the non-existence of mirror symmetry in some cases. In order to do this, we need the following previous result.

\begin{lemma} \label{inst-a}
Let $p(t),q(t)$ be two real polynomials, and let $A,B,C\in {\Bbb R}$. Also, let $R_{A,B,C}(t)=Ap(t)+Bq(t)+C$. Then there does not exist $(A_0,B_0)\in {\Bbb R}^2$ such that $R_{A_0,B_0,C}(t)$ has multiple roots for almost all $C\in {\Bbb R}$.
\end{lemma}

{\bf Proof.} Let $R_{A,B,C}'(t)$ denote the derivative of $R_{A,B,C}(t)$ with respect to $t$. For a particular value of $C\in {\Bbb R}$, $R_{A_0,B_0,C}(t)$ has multiple roots iff $R_{A_0,B_0,C}(t)=A_0p(t)+B_0q(t)+C$ and $R_{A_0,B_0,C}'(t)=A_0p'(t)+B_0q'(t)$ have a common factor. Moreover, $R_{A_0,B_0,C}(t)$ has multiple roots for almost all $C\in {\Bbb R}$ iff the former polynomials have a common factor regardless of the value of $C$. But this cannot happen because the second polynomial does not depend on $C$.
 \qed

Then we have the following result.

\begin{theorem} \label{th-prohib}
Let $r=\deg_t(x(t))$, $s=\deg_t(y(t))$. The following statements hold:
\begin{itemize}
\item [(1)] If $r$ (resp. $s$) is odd, then ${\mathcal C}$ cannot have any symmetry axis parallel to the $x$-axis (resp. $y$-axis).
\item [(2)] If $r=s$ and $r,s$ are odd, then ${\mathcal C}$ has at most one symmetry axis $Ax+By+C=0$ parallel to the vector $(a_r,a_s)$, where $a_r,a_s$ are the leading coefficients of $x(t),y(t)$, respectively.
    \end{itemize}
    \end{theorem}

    {\bf Proof.} Let us see (1). If ${\mathcal C}$ has a symmetry axis parallel to the $x$-axis (resp. the $y$-axis), then the family of lines $x=a$ (resp. $y=b$) must intersect ${\mathcal C}$ at an even number of points, for almost all $a\in {\Bbb R}$ (resp. $b\in {\Bbb R}$); however, if $r$ (resp. $s$) is odd this cannot happen. Now in order to see (2), one observes that if ${\mathcal C}$ has a symmetry axis ${\mathcal L}$, then the family of normal lines to ${\mathcal L}$ has the form $A_0x+B_0y+C=0$ for a certain $(A_0,B_0)\in {\Bbb R}^2$, where $C$ is a real parameter. The intersection of this family with ${\mathcal C}$ amounts to $Ax(t)+By(t)+C=0$, and this equation must have an even number of solutions for almost all $C\in {\Bbb R}$. Now assume that the degree of this equation is odd, and let $\delta$ stand for the number of real solutions of the equation; in that case, either $A_0a_r+B_0a_s=0$, or the equation has multiple solutions. However, by Lemma \ref{inst-a} the latter cannot happen for almost all values of $C$, and (2) follows. \qed

Example 2 in Subsection 2 provides an example for the statement (2) of this theorem. Now by
putting together Theorem \ref{conclu-some-pro} and Theorem \ref{th-prohib}, we get the following corollary of Theorem \ref{th-prohib}.

\begin{corollary} \label{no-symm}
If both $r=\deg_t(x(t))$ and $s=\deg_t(y(t))$ are odd, and $r\neq s$, then ${\mathcal C}$ does not exhibit mirror symmetry.
\end{corollary}

The results in Theorem \ref{conclu-some-pro}, Theorem \ref{th-prohib} and Corollary \ref{no-symm} are displayed in the table below, where we follow the notation of Theorem \ref{th-prohib}.

\begin{center}
\begin{tabular}{|c|l|}
\hline
${\bf r}$ odd, ${\bf s}$ odd & \begin{tabular}{c|l} ${\bf r<s:}$ & No mirror symmetry\\   ${\bf r=s:}$ & Axis parallel to $(a_r,a_s)$, if any. \\ ${\bf r>s}$ &  No mirror symmetry \end{tabular}\\ \hline
${\bf r}$ odd, ${\bf s}$ even & Axis parallel to $y$-axis, if any.\\ \hline
${\bf r}$ even, ${\bf s}$ odd & Axis parallel to $x$-axis, if any.\\ \hline
${\bf r}$ even, ${\bf s}$ even & \begin{tabular}{c|l} ${\bf r<s:}$ & Axis parallel to $y$-axis, if any.\\  ${\bf r=s:}$ & (nothing to say) \\  ${\bf r>s:}$ & Axis parallel to $x$-axis, if any. \end{tabular} \\ \hline \end{tabular}
\end{center}

        Finally, one might wonder whether central symmetry and mirror symmetry can coincide at the same time. Let us see that the answer, in our case, is negative. For this purpose, we first need the following result.

        \begin{proposition} \label{th-coincidence}
        Let ${\mathcal C}$ be an algebraic curve, and assume that it presents rotation symmetry, with center $P_0$ and angle $\theta$. If ${\mathcal L}$ is an axis of symmetry of ${\mathcal C}$, then $P_0\in {\mathcal L}$.
        \end{proposition}

        {\bf Proof.} Assume by contradiction that $P_0\notin {\mathcal L}$, and let $P_0'\neq P_0$ be the symmetric of $P_0$ w.r.t. ${\mathcal L}$. Furthermore, let $\tilde{\mathcal C}$ be the symmetric of ${\mathcal C}$ w.r.t. ${\mathcal L}$. Since ${\mathcal C}$ is by hypothesis symmetric w.r.t. ${\mathcal L}$, we have $\tilde{\mathcal C}={\mathcal C}$. Hence, a rotation of ${\mathcal C}$ around $P_0$ with angle $\theta$ amounts to a rotation of $\tilde{\mathcal C}$ around $P_0'$ with angle $-\theta$. But since $\tilde{\mathcal C}={\mathcal C}$, it holds that ${\mathcal C}$ has two symmetry centers, namely $P_0$ and $P_0'$. However this cannot happen because ${\mathcal C}$ is algebraic, and therefore it can have just one symmetry center (see Theorem 5.3 in \cite{LRTh}).
        So, $P_0\in {\mathcal L}$. \qed

        Now we can state the result.

        \begin{theorem} \label{no-both-sym}
       Polynomial curves cannot exhibit central symmetry and mirror symmetry at the same time.
       \end{theorem}

       {\bf Proof.} If ${\mathcal C}$ has a center of symmetry $P_0$ and a symmetry axis ${\mathcal L}$, then by Proposition \ref{th-coincidence} it holds that
      $P_0\in {\mathcal L}$. In that case, it is easy to see that ${\mathcal C}$ must also have another symmetry axis ${\mathcal L}'$, namely the normal line to ${\mathcal L}$ at $P_0$. But by Lemma \ref{one-most}, this is not possible. \qed

      \section{Algorithms Performance} \label{float}

In this section, we provide evidence of the practical performance of the algorithms derived in the preceding sections. We have implemented our algorithms in the computer algebra system Maple 14, and we have run several examples on a computer equipped with a 2GHz microprocessor, and 8Gb of RAM. The parametrizations of the examples tested can be checked in Appendix I. One may observe from the timings given that we can analyze curves of high degrees in seconds. Now in each case we provide the degree (i.e. the highest power of $t$ appearing in the components $x,y$ of the parametrization), the infinity norm of the coefficients in $x(t),y(t)$ (although we have always tried cases where the temptative value of $\beta$ was of reasonable size), and the timing in seconds. It is worth mentioning that the implicit equations of the considered curves are very costly to compute, and even impossible in some of the cases. Sometimes we
can get a quicker response by applying certain results (for instance, Theorem \ref{spec-case}, Theorem \ref{spec-case-mirror}, etc.) When that happens, we include an observation in the last column; some other details are also spelt there. Now the following table corresponds to examples where central symmetry was tested.

\hspace*{-0.5in}\begin{tabular}{|c|c|c|c|c|c|}
      \hline
      Ex. & Deg. & Norm & Cent. Sym. & Timing & Comments \\
      \hline
      2 & 83 & $2^{81}$ & Yes & $0.718$ &   \\
      3 & 3 & 15 & Yes & $0.343$ & Cubic with $c_2=0$ (see Th. \ref{spec-case})\\
      4 & 7 & 15  & Yes & $0.437$ & Curve of degree $\geq 4$ with $c_{n-1}=0$ (see Th. \ref{spec-case})\\
      37 & 21 & $2^{30}$ & No & $0.468$ & $\beta\notin {\Bbb R}$: with this, $0.172$ secs.\\
      38 & 45 & 194 & No & $0.421$ & $\beta\notin {\Bbb R}$: with this, $0.109$ secs.\\
      39 & 35 & $8.06\cdot 10^{9}$ & Yes &  $1.248$ & Algebraic coefficients. \\
      40 & 77 & 83 & No & $0.500$ & $\beta=0$ but $(-1,0)$ is not a solution of the system.\\
      41 & 95 & $5.89\cdot 10^{27}$ & Yes & $5.850$ & Algebraic coefficients. \\
      \hline
      \end{tabular}

Now we consider mirror symmetry.

      \hspace*{-0.6in}\begin{tabular}{|c|c|c|c|c|c|}
      \hline
      Ex. & Degree & Norm & Mirr. Sym. & Timing & Comments \\
      \hline
      6 & 21 & $2^{31}$ & Yes & $1.388$ &   \\
      7 & 21 & 3.5 & Yes & $0.983$ & Algebraic coefficients. \\
      8 & 69 & 276  & No & $3.229$ & Case $c_{n-1}=0$ (see Th. \ref{spec-case-mirror}; with this, $0.203$ secs.)\\
      9 & 45 & 461 & No & $1.653$ &  Case $c_{n-1}=0$ (see Th. \ref{spec-case-mirror}; with this, $0.203$ secs.)\\
      10 & 45 & 194 & No & $1.653$ & \\
      11 & 60 & 844 & Yes & $2.465$ & $r$ even, $s$ odd (see Th. \ref{conclu-some-pro})\\
      12 & 91 & $2^{89}$ & Yes & $15.803$ &  \\
      42 & 77 & $2^{75}$ & No & 9.219 &  $Q^{\star}(\beta)\notin {\Bbb R}$; with this, $9.017$ secs.\\
      43 & 35& $2^{38}$& No& $2.901$ & $r,s$ odd (see Th. \ref{conclu-some-pro}; with this, 0 secs.) \\
      44 & 56& $2^{52}$& Yes& $23.385$ & Algebraic coefficients. \\
      \hline
      \end{tabular}

The complexity analysis of the above algorithms is provided in the following theorem.

\begin{theorem} \label{complexity}
The complexity of the above algorithms for detecting central and mirror symmetry is ${\mathcal O}(n^3)$, where $n$ is the degree of the parametrization, $\varphi(t)$.
\end{theorem}

{\bf Proof.} We prove the result for the case of central symmetry. The mirror symmetry case is analogous. Now the complexity is dominated by that of evaluating the equations $(1), (2), \ldots, (n)$ of the system ${\mathcal S}$ (see equation \ref{eq-central} in Subsection \ref{subsec-detect}) at $\alpha=-1,\beta=\frac{-2c_{n-1}}{n\cdot c_n}$. The equation (i) of ${\mathcal S}$ has the form
\begin{equation}\label{compl}
-c_i=\alpha^i\cdot \sum_{k=i}^n c_k {k \choose i} \beta^{k-i}
\end{equation}
for $i=1,\ldots,n$. So, it consists of the sum of $n-i+2$ terms: at the left hand-side of (\ref{compl}) we have just one term; at the right hand-side we have $n-i+1$ terms, and the $j$-th term (where $j=0,\ldots,n-i$) is the product of $j+3$ complex numbers. Since
\[\sum_{j=0}^{n-i}(j+3)=\frac{n-i+6}{2}(n-i+1)={\mathcal O}(n^2),\]
 in order to evaluate the equation (i) for $\alpha=-1,\beta=\frac{-2c_{n-1}}{n\cdot c_n}$ we need to perform ${\mathcal O}(n^2)$ multiplications, and ${\mathcal O}(n)$ sums. We have to do this for $n$ equations; so, we get a complexity ${\mathcal O}(n^3)$. \qed

\section{Conclusions and Further Work} \label{conclusion}

Here we have presented an efficient method for detecting the symmetries, and computing them in the affirmative case, of algebraic curves defined by means of polynomial parametrizations. The method is a combination of geometric and algebraic ideas: on the one hand, we observe that the nature of the symmetry leads to a new parametrization of the curve; on the other hand, whenever we start from a proper parametrization, this second parametrization must also be proper, and is related with the first one by means of a linear mapping (see Lemma \ref{z-prop}). By using this and writing our parametrization in complex form (which is really useful here) we obtain a polynomial system with three (complex) unknowns but with a triangular structure, that can therefore be solved in a fast and efficient way. In fact, we provide closed expressions for the temptative solutions of the system, and for the elements defining the symmetries of ${\mathcal C}$, if any. The complexity analysis and the experimentation performed so far prove evidence of the efficiency of the method. Now it is natural to wonder if this idea can be extended for polynomial curves in higher dimensions. As for central symmetry, this is certainly the case: even if we cannot use complex numbers anymore, we can always split the problem into the analysis of several, triangular, polynomial systems with a similar behavior to the one observed here. But the situation is more complicated when it comes to other symmetry types: in ${\Bbb R}^3$, for instance, we may have symmetry with respect to a plane, or symmetry with respect to an axis (the case when the curve is invariant with respect to rotations about a line); in both cases, the generalization is not direct, and is in fact intended to be addressed in a future work.


\newpage
\section{Appendix I} \label{Appendix}

For each curve, we list: the number assigned to the example in our database, the parametrization of the curve, and some observations.

{\bf Example 2:}
\[
\left\{\begin{array}{l}
x(t)=3/2+3(t+1)^{83}-(t+1)^{53}+3(t+1)^{11}-(t+1)^{63}+t\\
y(t)=2/3+(t+1)^{83}+(t+1)^{53}-3(t+1)^{11}-(t+1)^{63}+t
\end{array}\right.
\]

\begin{itemize}
\item Degree of the parametrization: $83$
\item Symmetry tested: central.
\item Values for $\alpha,\beta$: $\alpha=-1,\beta=-2$.
\end{itemize}

{\bf Example 3:}
\[
\left\{\begin{array}{l}
x(t)=-6+15t^3-4t\\
y(t)=-2+2t^3-t
\end{array}\right.
\]
\begin{itemize}
\item Degree of the parametrization: $3$
\item Symmetry tested: central.
\item Values for $\alpha,\beta$: $\alpha=-1,\beta=0$.
\item Observations: we can get the answer in 0 seconds by applying Theorem \ref{spec-case}.
\end{itemize}

{\bf Example 4:}
\[
\left\{\begin{array}{l}
x(t)=15t^7-4t^5-t-2\\
y(t)=2t^7-t^5-t^3+t+1
\end{array}\right.
\]
\begin{itemize}
\item Degree of the parametrization: $7$
\item Symmetry tested: central.
\item Values for $\alpha,\beta$: $\alpha=-1,\beta=0$.
\item Observations: we can get the answer in 0 seconds by applying Theorem \ref{spec-case}.
\end{itemize}

{\bf Example 6:}
\[
\left\{\begin{array}{l}
x(t)=(2t+1)^{20}+(2t+1)^{18}+(2t+1)^{10}+1+(2t+1)^{21}-3(2t+1)^5+(2t+1)^3\\
y(t)=-(2t+1)^{20}-(2t+1)^{18}-(2t+1)^{10}-1+(2t+1)^{21}-3(2t+1)^5+(2t+1)^3
\end{array}\right.
\]
\begin{itemize}
\item Degree of the parametrization: $21$
\item Symmetry tested: mirror.
\item Values for $\alpha,\beta$: $\alpha=-1,\beta=-1$.
\end{itemize}

{\bf Example 7:}

Let
\[
\begin{array}{l}
u(t)=t^{21}-3t^5+t^3\\
v(t)=t^{20}+t^{18}+t^{10}+1
\end{array}
\]
Then the curve is
\[
\left\{\begin{array}{l}
x(t)=5+\frac{\sqrt{3}}{2}u(t)-\frac{1}{2}v(t)\\
y(t)=-3+\frac{1}{2}u(t)+\frac{\sqrt{3}}{2}v(t)
\end{array}\right.
\]
\begin{itemize}
\item Degree of the parametrization: $21$
\item Symmetry tested: mirror.
\item Values for $\alpha,\beta$: $\alpha=-1,\beta=0$.
\item Observations: even though the parametrization contains algebraic coefficients, the algorithm works perfectly.
\end{itemize}

{\bf Example 8:}
\[
\hspace{-1in}\left\{\begin{array}{rcl}
x(t)&=&-219t^{69}+100t^{55}-46t^{44}-12t^{37}-249t^{31}-30t^{30}-150t^{26}+184t^{17}+186t^{16}-12t^{13}-\\
&&-148t^9-246t^6+5\\
y(t)&=&-146t^{69}-150t^{55}+69t^{44}-8t^{37}-166t^{31}-20t^{30}+225t^{26}-276t^{17}+124t^{16}+18t^{13}+\\
&&222t^9-164t^6-3
\end{array}\right.
\]
\begin{itemize}
\item Degree of the parametrization: $69$
\item Symmetry tested: mirror.
\item Values for $\alpha,\beta$: $\alpha=-1,\beta=0$.
\item Observations: by using Theorem \ref{spec-case-mirror} we can get a very fast answer. 
\end{itemize}

{\bf Example 9:}
\[
\hspace{-1in}\left\{\begin{array}{l}
x(t)=348t^{45}+246t^{34}-240t^{25}-224t^{22}+132t^{19}-461t^{17}+51t^{14}+225t^{11}+388t^6-292t^4+5\\
y(t)=261t^{45}-328t^{34}+320t^{25}-168t^{22}-176t^{19}+98t^{17}-68t^{14}-300t^{11}+291t^6-219t^4-3
\end{array}\right.
\]
\begin{itemize}
\item Degree of the parametrization: $45$
\item Symmetry tested: mirror.
\item Values for $\alpha,\beta$: $\alpha=-1,\beta=0$.
\item Observations: by using Theorem \ref{spec-case-mirror} we can get a very fast answer.
\end{itemize}

{\bf Example 10:}
\[
\hspace{-1in}\left\{\begin{array}{l}
x(t)=-t^{44}+87t^{25}-56t^{18}-62t^8+97t^5-73t^2+5\\
y(t)=174t^{45}+t^{44}-87t^{25}-112t^{22}+56t^{18}-124t^{17}+62t^8+194t^6-97t^5-146t^4+73t^2-3
\end{array}\right.
\]
\begin{itemize}
\item Degree of the parametrization: $45$
\item Symmetry tested: mirror.
\item Values for $\alpha,\beta$: $\alpha=-1,\beta=-\frac{1}{3915}$.
\end{itemize}

{\bf Example 11:}
\[
\left\{\begin{array}{l}
x(t)=1+t^{60}-94t^{30}-56t^{14}-62t^4\\
y(t)=2-t^{59}-82t^{33}+844t^{27}-17t^{19}
\end{array}\right.
\]
\begin{itemize}
\item Degree of the parametrization: $60$
\item Symmetry tested: mirror.
\item Values for $\alpha,\beta$: $\alpha=-1,\beta=0$.
\item Observation: corresponds to the case $r$ even, $s$ odd (see Theorem \ref{conclu-some-pro})
\end{itemize}

{\bf Example 12:}
\[
\hspace{-1in}\left\{\begin{array}{rcl}
x(t)&=&-2+2(t+1)^{90}+(t+1)^{88}+73(t+1)^{86}-4(t+1)^{54}-83(t+1)^{48}-82(t+1)^{20}+(t+1)^{91}+\\
&&+82(t+1)^{73}-12(t+1)^{47}+17(t+1)^{29}\\
y(t)&=&4+2(t+1)^{90}+(t+1)^{88}+73(t+1)^{86}-4(t+1)^{54}-83(t+1)^{48}-82(t+1)^{20}-(t+1)^{91}-\\
&&-82(t+1)^{73}+12(t+1)^{47}-17(t+1)^{29}
\end{array}\right.
\]
\begin{itemize}
\item Degree of the parametrization: $91$
\item Symmetry tested: mirror.
\item Values for $\alpha,\beta$: $\alpha=-1,\beta=-2$.
\end{itemize}

{\bf Example 37:}
\[
\left\{\begin{array}{l}
x(t)=1+(2t+1)^{20}+(2t+1)^{18}+(2t+1)^{10}+(2t+1)^{21}-3(2t+1)^5+(2t+1)^3\\
y(t)=-1-(2t+1)^{20}-(2t+1)^{18}-(2t+1)^{10}+(2t+1)^{21}-3(2t+1)^5+(2t+1)^3
\end{array}\right.
\]
\begin{itemize}
\item Degree of the parametrization: $20$
\item Symmetry tested: central.
\item Values for $\alpha,\beta$: $\alpha=-1,\beta=-1+\frac{1}{21}i$.
\item Observations: Since $\beta\notin {\Bbb R}$, we conclude rapidly that the curve has not central symmetry.
\end{itemize}

{\bf Example 38:}
\[
\hspace{-0.7in}\left\{\begin{array}{l}
x(t)=-t^{44}+87t^{25}-56t^{18}-62t^8+97t^5-73t^2+5\\
y(t)=174t^{45}+t^{44}-87t^{25}-112t^{22}+56t^{18}-124t^{17}+62t^8+194t^6-97t^5-146t^4+73t^2-3
\end{array}\right.
\]
\begin{itemize}
\item Degree of the parametrization: $45$
\item Symmetry tested: central.
\item Values for $\alpha,\beta$: $\alpha=-1,\beta=-\frac{1}{3915}-\frac{1}{3915}i$.
\item Observations: Since $\beta\notin {\Bbb R}$, we conclude rapidly that the curve has not central symmetry.
\end{itemize}

{\bf Example 39:}

Let
\[
\begin{array}{l}
u(t)=2(t+1)^{35}-3(t+1)^{31}+2(t+1)^{27}-5(t+1)^{23}+3(t+1)^{15}-(t+1)\\
v(t)=2(t+1)^{33}-3(t+1)^{31}+2(t+1)^{27}-5(t+1)^{19}+3(t+1)^{15}-(t+1)^7
\end{array} 
\]
Then the curve is 
\[
\left\{\begin{array}{l}
x(t)=\frac{\sqrt{2}}{2}+\frac{\sqrt{3}}{2}u(t)+\frac{1}{2}v(t)\\
y(t)=-\frac{\sqrt{2}}{2}-\frac{1}{2}u(t)+\frac{\sqrt{3}}{2}v(t)
\end{array}\right.
\]
\begin{itemize}
\item Degree of the parametrization: $35$
\item Symmetry tested: central.
\item Values for $\alpha,\beta$: $\alpha=-1,\beta=-2$.
\item Observations: even though the parametrization contains algebraic coefficients, the algorithm works perfectly.
\end{itemize}

{\bf Example 40:}
\[
\left\{\begin{array}{l}
x(t)=4t^{77}-83t^{52}+t^{28}-62t^{23}-30t^{15}+t\\
y(t)=2t^{63}-41t^{32}+3t^{27}-28t^{19}+4t^{15}-t^7
\end{array}\right.
\]
\begin{itemize}
\item Degree of the parametrization: $77$
\item Symmetry tested: central.
\item Values for $\alpha,\beta$: $\alpha=-1,\beta=0$.
\item Observations: even though $\beta$ is real, the pair $(-1,0)$ is not a solution of the system. So, the curve has not central symmetry. 
\end{itemize}

{\bf Example 41:}

Let
\[
\begin{array}{l}
u(t)=2(t+1)^{95}-3(t+1)^{91}+2(t+1)^{87}-5(t+1)^{23}+3(t+1)^{15}-(t+1)\\
v(t)=2(t+1)^{93}-3(t+1)^{71}+2(t+1)^{57}-5(t+1)^{39}+3(t+1)^{15}-(t+1)^7
\end{array}
\]
Then the curve is
\[
\left\{\begin{array}{l}
x(t)=\frac{\sqrt{2}}{2}+\frac{\sqrt{3}}{2}u(t)+\frac{1}{2}v(t)\\
y(t)=-\frac{\sqrt{2}}{2}-\frac{1}{2}u(t)+\frac{\sqrt{3}}{2}v(t)
\end{array}\right.
\]
\begin{itemize}
\item Degree of the parametrization: $95$
\item Symmetry tested: central.
\item Values for $\alpha,\beta$: $\alpha=-1,\beta=-2$.
\item Observations: even though the parametrization contains algebraic coefficients, the algorithm works perfectly.
\end{itemize}

{\bf Example 42:}
\[
\hspace{-1.2in}\left\{\begin{array}{l}
x(t)=3(t+1)^{77}+(t+1)^{76}-73(t+1)^{71}-4(t+1)^{39}-83(t+1)^{33}-10(t+1)^{32}+62(t+1)^{18}-82(t+1)^8
\\
y(t)=(t+1)^{77}-(t+1)^{76}-73(t+1)^{70}-4(t+1)^{38}-83(t+1)^{32}-10(t+1)^{31}+62(t+1)^{17}-82(t+1)^7
\end{array}\right.
\]
\begin{itemize}
\item Degree of the parametrization: $77$
\item Symmetry tested: mirror.
\item Values for $\alpha,\beta$: $\alpha=-1,\beta=-\frac{772}{385}$.
\item Observations: one may see that $Q^{\star}(\beta)\notin {\Bbb R}$. 
\end{itemize}

{\bf Example 43:}
\[
\left\{\begin{array}{l}
x(t)=-82-73(t+1)^{31}-4(t+1)^{30}-83(t+1)^7-10(t+1)^4+62(t+1)^2
\\
y(t)=95(t+1)^{35}+11(t+1)^{30}-49(t+1)^{29}-47(t+1)^{23}+40(t+1)^{20}-81(t+1)^8
\end{array}\right.
\]
\begin{itemize}
\item Degree of the parametrization: $35$
\item Symmetry tested: mirror.
\item Values for $\alpha,\beta$: $\alpha=-1,\beta=-2$.
\item Observations: one may see that $r,s$ are odd (see Theorem \ref{conclu-some-pro}). By using this, we get an answer in 0 seconds.
\end{itemize}

{\bf Example 44:}

Let
\[
\begin{array}{l}
u(t)= t^{56}-t^{40}+2t^{28}-3t^{12}+t^4-t^2\\
v(t)=t^{53}-t^{41}+2t^{23}-3t^{11}+t^3-t^5
\end{array}
\]
Then the curve is
\[
\left\{\begin{array}{l}
x(t)=\sqrt{2}+\frac{\sqrt{3}{2}}u(t)+\frac{1}{2}v(t)\\
y(t)=\sqrt{3}-\frac{1}{2}u(t)+\frac{\sqrt{3}}{2}v(t)
\end{array}\right.
\]
\begin{itemize}
\item Degree of the parametrization: $56$
\item Symmetry tested: mirror.
\item Values for $\alpha,\beta$: $\alpha=-1,\beta=-2$.
\item Observations: even though the parametrization contains algebraic coefficients, the algorithm works perfectly.
\end{itemize}

\end{document}